\documentclass[12pt]{article}
\usepackage{amsfonts,amssymb}

\def\C{\mathbb{C}}
\def\N{\mathbb{N}}
\def\Z{\mathbb{Z}}

\def\bq{ \begin{equation} }
\def\eq{ \end{equation} }
\def\ben{ \begin{eqnarray} }
\def\en{ \end{eqnarray} }

\def\frac#1#2{{#1\over #2}}

\def\on#1#2{\mathop{\vbox{\ialign{##\crcr\noalign{\kern2pt}
$\scriptstyle{#2}$\crcr\noalign{\kern2pt\nointerlineskip}
\kern-2pt$\hfil\displaystyle{#1}\hfil$\crcr}}}\limits}

\begin{document}

\baselineskip=15pt
\vspace{1cm} \noindent {\LARGE \textbf {A family of
 (2+1)-dimensional hydrodynamic type systems possessing  pseudopotential }}
\vskip1cm \hfill
\begin{minipage}{13.5cm}
\baselineskip=15pt {\bf Alexander Odesskii }

 {\footnotesize
Landau Institute for Theoretical Physics, Moscow, Russia
\\
  School of Mathematics, The University of
Manchester, UK}\\
\vskip1cm{\bf Abstract}

We construct a family of integrable hydrodynamic type systems with
three independent and $n\ge 2$ dependent variables in terms of
solutions of linear system of PDEs with rational coefficients. We
choose the existence of a pseudopotential as a criterion of
integrability. In the case $n=2$ this family is a general solution
of the classification problem for such systems. We give also an
elliptic analog of this family in the case $n>2$.

\end{minipage}

\vglue1cm \textbf{Address}: Landau Institute for Theoretical
Physics, Kosygina 2, 119334, Moscow, Russia

\textbf{E-mail}: alexander.odesskii@manchester.ac.uk
\newpage

\centerline{\Large\bf Introduction}
\medskip

We address the problem of classification of integrable
$(2+1)$-dimensional quasilinear systems
 \begin{equation}   \label{fergener}
 {\bf u}_t=A({\bf u})\, {\bf u}_x+ B({\bf u}) \, {\bf u}_y,
  \end{equation}
where $t,x,y$ are independent variables, ${\bf u}$ is an
$n$-component column vector, and $A({\bf u})$ and $B({\bf u})$ are
$n\times n$-matrices. A general theory of such systems was developed
in the papers \cite{ferhus1, ferhus2, ferhus3}. This theory is based
on the existence of sufficiently many of the hydrodynamic reductions
\cite{gibtsar, ferhus1} which has been proposed as the definition of
integrability. In the first nontrivial case $n=2$ the complete set
of integrability conditions has been found in the paper
\cite{ferhus2} in the form of a complicated system of PDEs for the
entries of the matrices $A$ and $B$. In the case $n>2$ the complete
set of integrability conditions in terms of the matrices $A$ and $B$
is not known. However, all known examples of integrable
(2+1)-dimensional hydrodynamic type systems possess the so-called
scalar pseudopotential\footnote{This means that the overdetermined
system (\ref{pseudo0}) for $\Psi$ is compatible if and only if ${\bf
u}$ is a solution of (\ref{fergener}). }
\begin{equation}\label{pseudo0}
 \Psi_t=f(\Psi_y, {\bf u}), \qquad   \Psi_x=g(\Psi_y, {\bf u}).
\end{equation}
Moreover, it was proven in \cite{ferhus2} that for $n=2$ the
integrability conditions are equivalent to the existence of the
scalar pseudopotential. The scalar pseudopotential plays an
important role in the theory of the universal Whitham hierarchy
\cite{kr1,kr3,kr4}. Various integrable systems possessing
pseudopotential were constructed and studied in \cite{zakh},
\cite{dun} and many other papers. In the paper \cite{odsok} the case
of a constant matrix $A$ was considered and hydrodynamic type
systems (\ref{fergener}) possessing a pseudopotential with movable
singularities were classified.

In this paper we give the general solution to the classification
problem for the integrable hydrodynamic type systems
(\ref{fergener}) in the case $n=2$. It turns out that the answer can
be written in terms of three arbitrary linear independent solutions
of a certain system of linear PDEs with rational coefficients and a
3-dimensional space of solutions. The action of the group $GL_3$ on
this space corresponds to the action of the same group on the space
spanned by the independent variables $t,x,y$. It is also possible to
construct a similar family of hydrodynamic type systems possessing a
pseudopotential for higher $n.$ The system is written in terms of
three arbitrary linear independent solutions of a certain system of
linear PDEs with rational coefficients and an $n+1$-dimensional
space of solutions.

Let us describe the contents of the paper. In Section 1 we recall
the background material. In Section 2 we construct a family of
hydrodynamic type systems possessing a pseudopotential for arbitrary
$n\geq2$. In Section 3 we recall the results of the paper
\cite{ferhus2} on two components integrable systems and prove that
our family in the case $n=2$ is a general solution for the
classification problem. In Section 4 we give an elliptic analog of
our construction.

\section{Hydrodynamic type systems and their pseudopotentials}

Consider a $(2+1)$-dimensional quasilinear system
\begin{equation}   \label{gener}
u_{it}=\sum_{j=1}^n a_{ij}({\bf u})\,u_{j x}+ \sum_{j=1}^n
b_{ij}({\bf u})\,u_{j y}, ~~~i=1,...,n.
\end{equation}
This system is called integrable \cite{ferhus1} if it possesses
`sufficiently many' exact solutions of the form ${\bf u}={\bf
u}(R^1, ..., R^N)$ where the so-called Riemann invariants $R^1, ...,
R^N$ solve a pair of commuting diagonal systems
\begin{equation}
R^i_t=\lambda^i(R)\ R^i_y, ~~~~ R^i_x=\mu^i(R)\ R^i_y, \label{R}
\end{equation}
and the number $N$ of Riemann invariants is allowed to be arbitrary.

A pair of equations of the form
\begin{equation}\label{pseudo}
 \Psi_t=f(\Psi_y, u_1,\dots, u_n), \qquad   \Psi_x=g(\Psi_y, u_1,\dots, u_n),
\end{equation}
with respect to unknown $\Psi$ is called a pseudopotential for
equation (\ref{gener}) if the compatibility  condition
$\Psi_{tx}=\Psi_{xt}$ for (\ref{pseudo}) is equivalent to
(\ref{gener}). Differentiating (\ref{pseudo}), we find that this
compatibility condition is given by
\begin{equation}\label{comp}
f_{\xi}\, \sum_{i=1}^n u_{iy} \, g_{u_i}   + \sum_{i=1}^n u_{ix}\,
f_{u_i} = g_{\xi} \sum_{i=1}^n u_{iy} \, f_{u_i} + \sum_{i=1}^n
u_{it} \, g_{u_i}.
\end{equation}
Here and below we denote $\Psi_y$ by $\xi$. Substituting the right
hand side of (\ref{gener}) for $t$-derivatives and splitting with
respect to $x$ and $y$-derivatives, we get that for any $i$ the
following relations hold:
\begin{equation}\label{psecon1}
 f_{u_i}=\sum_{j=1}^n a_{ji}  g_{u_j},
 \end{equation}
 \begin{equation}\label{psecon2}     f_{\xi} \,
 g_{u_i} - g_{\xi} \,  f_{u_i}=\sum_{j=1}^n b_{ji}  g_{u_j}.
\end{equation}

We will use the following

{\bf Lemma 1.} A pair of equations (\ref{pseudo}) is a
pseudopotential for some system of the form (\ref{gener}) iff the
functions $\{f_{u_i},g_{u_i},f_{\xi} \,
 g_{u_i} - g_{\xi}  f_{u_i};i=1,...,n\}$ constitute an
 $n$-dimensional linear subspace in the space of functions in $\xi$
 and $\{g_{u_i};i=1,...,n\}$ is a basis in this space.

 {\bf Proof.} Indeed, $\{g_{u_i};i=1,...,n\}$ is a basis in this
 space of functions in $\xi$ iff $f_{u_i},f_{\xi} \,
 g_{u_i} - g_{\xi}  f_{u_i}$ can be uniquely written in the form
 (\ref{psecon1}), (\ref{psecon2}) where coefficients $a_{ji}$,
 $b_{ji}$ does not depend on $\xi$.

 In some cases it is more convenient to define a pseudopotential in
 parametric form
\begin{equation}\label{pseudopar}
 \Psi_y=F(\zeta, u_1,\dots, u_n), \qquad   \Psi_t=G(\zeta, u_1,\dots, u_n), \qquad   \Psi_x=H(\zeta, u_1,\dots, u_n)
\end{equation}
which takes the form (\ref{pseudo}) if one expresses the parameter
$\zeta$ in terms of $\Psi_y$ from the first equation and substitutes
into the second and the third equations. The form (\ref{pseudopar})
is more symmetric with respect to the variables $t,x,y$. In
particular, linear change of the variables $y,t,x$ corresponds to
the same linear change of the functions $F,G,H$. The compatibility
conditions for the system (\ref{pseudopar}) read
\begin{equation}\label{comppar}
\sum_{i=1}^n(H_{\zeta}G_{u_i}-G_{\zeta}H_{u_i})u_{iy}+\sum_{i=1}^n(F_{\zeta}H_{u_i}-H_{\zeta}F_{u_i})u_{it}+
\sum_{i=1}^n(G_{\zeta}F_{u_i}-F_{\zeta}G_{u_i})u_{ix}=0.
\end{equation}
We can slightly generalize Lemma 1 in the following way:

{\bf Lemma 2.} The compatibility conditions for the system
(\ref{pseudopar}) are equivalent to a quasilinear system of $m$
linear independent equations the form
\begin{equation}   \label{genern}
\sum_{j=1}^n a_{ij}({\bf u})\,u_{jy}+\sum_{j=1}^n b_{ij}({\bf
u})\,u_{j t}+ \sum_{j=1}^n c_{ij}({\bf u})\,u_{j x}=0, ~~~i=1,...,m
\end{equation}
iff the functions
$\{H_{\zeta}G_{u_i}-G_{\zeta}H_{u_i},
F_{\zeta}H_{u_i}-H_{\zeta}F_{u_i},
G_{\zeta}F_{u_i}-F_{\zeta}G_{u_i}; i=1,...,n\}$ constitute an
$m$-dimensional linear subspace in the space of functions in
$\zeta$.

{\bf Proof.} Let $\{S_1(\zeta),...,S_m(\zeta)\}$ be a basis in this
linear space and
$$H_{\zeta}G_{u_i}-G_{\zeta}H_{u_i}=\sum_{j=1}^ma_{ji}({\bf u})S_j,~~~
F_{\zeta}H_{u_i}-H_{\zeta}F_{u_i}=\sum_{j=1}^mb_{ji}({\bf u})S_j,~~~
G_{\zeta}F_{u_i}-F_{\zeta}G_{u_i}=\sum_{j=1}^mc_{ji}({\bf u})S_j.$$
Substituting these equations into (\ref{comppar}) and equating to
zero coefficients at $S_1,...,S_m$ we obtain (\ref{genern}).

\section{Construction of a family of $n$-component systems possessing pseudopotential}

Define a function $F(\zeta, u_1,\dots, u_n)$ as a solution of the
following systems of PDEs
$$F_{\zeta}=\phi(\zeta)\cdot\zeta^{-s_1}(\zeta-1)^{-s_2}(\zeta-u_1)^{-s_3}...(\zeta-u_n)^{-s_{n+2}},$$
\begin{equation}   \label{pspar}
F_{u_i}=-\frac{\phi(u_i)}{\zeta-u_i}\cdot\frac{\zeta^{1-s_1}(\zeta-1)^{1-s_2}(\zeta-u_1)^{1-s_3}...(\zeta-u_n)^{1-s_{n+2}}}
{u_i(u_i-1)(u_i-u_1)...\hat{i}...(u_i-u_n)},~~~i=1,...,n.
\end{equation}
Here $s_1,...,s_{n+2}$ are constants and
$$\phi(\zeta)=\alpha_0(u_1,...,u_n)+\alpha_1(u_1,...,u_n)\zeta+...+\alpha_n(u_1,...,u_n)\zeta^n$$
is a polynomial of degree $n$. The notation $\hat{i}$ means that the
$i$th multiplier is omitted in the product. The system (\ref{pspar})
is in involution iff the polynomial $\phi$ satisfies the following
system of PDEs
\begin{equation}   \label{lin}
\phi_{u_i}(\zeta)=\phi(u_i)\frac{\zeta(\zeta-1)(\zeta-u_1)...\hat{i}...(\zeta-u_n)}{u_i(u_i-1)(u_i-u_1)...\hat{i}...(u_i-u_n)}\times
\end{equation}
$$\left(\frac{s_1-1}{\zeta}+\frac{s_2-1}{\zeta-1}
+\frac{s_3-1}{\zeta-u_1}+...+\frac{s_{i+2}}{\zeta-u_i}+...+\frac{s_{n+2}-1}{\zeta-u_n}\right)-\frac{s_{i+2}}{\zeta-u_i}\phi(\zeta),~~~
i=1,...,n.$$ Here
$\phi_{u_i}(\zeta)=\alpha_{0u_i}+\alpha_{1u_i}\zeta+...+\alpha_{nu_i}\zeta^n$.
It is clear that if $\phi$ is a polynomial of degree $n$, then the
right hand side of (\ref{lin}) is also a polynomial of degree $n$.
Therefore, (\ref{lin}) is a well-defined system of linear PDEs for
coefficients of $\phi$. It can be checked straightforwardly that
this system is in involution. Therefore, there are $n+1$ linear
independent solutions.

Let $\phi$, $\phi_1$ and $\phi_2$ be three linear independent
solutions of the system (\ref{lin}). We assume that $\phi$,
$\phi_1$, $\phi_2$ are polynomials of degree $n$ with respect to
$\zeta$. Define functions $G(\zeta, u_1,\dots, u_n)$ and $H(\zeta,
u_1,\dots, u_n)$ similarly to (\ref{pspar}) by
$$G_{\zeta}=\phi_1(\zeta)\cdot\zeta^{-s_1}(\zeta-1)^{-s_2}(\zeta-u_1)^{-s_3}...(\zeta-u_n)^{-s_{n+2}},$$
\begin{equation}   \label{pspar1}
G_{u_i}=-\frac{\phi_1(u_i)}{\zeta-u_i}\cdot\frac{\zeta^{1-s_1}(\zeta-1)^{1-s_2}(\zeta-u_1)^{1-s_3}...(\zeta-u_n)^{1-s_{n+2}}}
{u_i(u_i-1)(u_i-u_1)...\hat{i}...(u_i-u_n)},~~~i=1,...,n
\end{equation}
for the function $G$ and
$$H_{\zeta}=\phi_2(\zeta)\cdot\zeta^{-s_1}(\zeta-1)^{-s_2}(\zeta-u_1)^{-s_3}...(\zeta-u_n)^{-s_{n+2}},$$
\begin{equation}   \label{pspar2}
H_{u_i}=-\frac{\phi_2(u_i)}{\zeta-u_i}\cdot\frac{\zeta^{1-s_1}(\zeta-1)^{1-s_2}(\zeta-u_1)^{1-s_3}...(\zeta-u_n)^{1-s_{n+2}}}
{u_i(u_i-1)(u_i-u_1)...\hat{i}...(u_i-u_n)},~~~i=1,...,n
\end{equation}
for the function $H$.

{\bf Proposition 1.} If the functions $F,G,H$ are defined by
(\ref{pspar}), (\ref{pspar1}) and (\ref{pspar2}), then the system
(\ref{pseudopar}) defines a pseudopotential for some system of the
form (\ref{genern}) with $m=n$.

{\bf Proof.} Equations (\ref{pspar}), (\ref{pspar1}), (\ref{pspar2})
imply
$$H_{\zeta}G_{u_i}-G_{\zeta}H_{u_i}=\vartheta_i(\zeta)\cdot\zeta^{1-2s_1}(\zeta-1)^{1-2s_2}(\zeta-u_1)^{1-2s_3}...(\zeta-u_n)^{1-2s_{n+2}},$$
\begin{equation}   \label{exp1}
F_{\zeta}H_{u_i}-H_{\zeta}F_{u_i}=\nu_i(\zeta)\cdot\zeta^{1-2s_1}(\zeta-1)^{1-2s_2}(\zeta-u_1)^{1-2s_3}...(\zeta-u_n)^{1-2s_{n+2}},
\end{equation}
$$G_{\zeta}F_{u_i}-F_{\zeta}G_{u_i}=\mu_i(\zeta)\cdot\zeta^{1-2s_1}(\zeta-1)^{1-2s_2}(\zeta-u_1)^{1-2s_3}...(\zeta-u_n)^{1-2s_{n+2}}$$
where functions $\vartheta_i(\zeta),\nu_i(\zeta),\mu_i(\zeta)$ are
defined by
$$\vartheta_i(\zeta)=\frac{1}{u_i(u_i-1)(u_i-u_1)...\hat{i}...(u_i-u_n)}\cdot\frac{\phi_1(u_i)\phi_2(\zeta)-\phi_2(u_i)\phi_1(\zeta)}
{\zeta-u_i},$$
\begin{equation}   \label{mu}
\nu_i(\zeta)=\frac{1}{u_i(u_i-1)(u_i-u_1)...\hat{i}...(u_i-u_n)}\cdot\frac{\phi_2(u_i)\phi(\zeta)-\phi(u_i)\phi_2(\zeta)}{\zeta-u_i},
\end{equation}
$$\mu_i(\zeta)=\frac{1}{u_i(u_i-1)(u_i-u_1)...\hat{i}...(u_i-u_n)}\cdot\frac{\phi(u_i)\phi_1(\zeta)-\phi_1(u_i)\phi(\zeta)}{\zeta-u_i}.$$
It is clear that $\vartheta_i(\zeta)$ $\nu_i(\zeta)$, $\mu_i(\zeta)$
are polynomials in the variable $\zeta$ of degree $n-1$. Therefore,
the space of functions $\{H_{\zeta}G_{u_i}-G_{\zeta}H_{u_i},
F_{\zeta}H_{u_i}-H_{\zeta}F_{u_i},
G_{\zeta}F_{u_i}-F_{\zeta}G_{u_i}; i=1,...,n\}$ in the variable
$\zeta$ is isomorphic to the space of polynomials in $\zeta$ of
degree less or equal to $n-1$. This space is $n$-dimensional and we
can apply Lemma 2.

Let us construct the system possessing pseudopotential defined by
(\ref{pspar}), (\ref{pspar1}), (\ref{pspar2}) explicitly in terms of
three linear independent solutions of the linear system (\ref{lin}).

{\bf Proposition 2.} Let
$$\phi(\zeta)=\alpha_0+\alpha_1\zeta+...+\alpha_n\zeta^n,$$
$$\phi_1(\zeta)=\beta_0+\beta_1\zeta+...+\beta_n\zeta^n,$$
$$\phi_2(\zeta)=\gamma_0+\gamma_1\zeta+...+\gamma_n\zeta^n.$$ Then the
system with pseudopotential defined by (\ref{pspar}),
(\ref{pspar1}), (\ref{pspar2}) can be written in the form
\begin{equation}   \label{ex}
\sum_{i,j,k}\frac{u_i^{j+l}}{(u_i-u_1)...\hat{i}...(u_i-u_n)}((\gamma_j\alpha_k-\gamma_k\alpha_j)u_{it}+
(\alpha_j\beta_k-\alpha_k\beta_j)u_{ix}+(\beta_j\gamma_k-\beta_k\gamma_j)u_{iy})=0.
\end{equation}
Here summation is made by $i,j,k$ subject to the constrains $1\leq
i\leq n$, $0\leq j\leq l<k\leq n$ with fixed $l$. For each
$l=1,...,n$ we have an equation.

{\bf Proof.} Substituting (\ref{exp1}) into (\ref{comppar}) we
obtain
\begin{equation}   \label{comp1}
 \sum_{i=1}^n \nu_i(\zeta)u_{it}+\sum_{i=1}^n
\mu_i(\zeta)u_{ix}+\sum_{i=1}^n \vartheta_i(\zeta)u_{iy}=0.
\end{equation}
Calculating polynomials
$\nu_i(\zeta),\mu_i(\zeta),\vartheta_i(\zeta)$ in terms of
coefficients of polynomials $\phi(\zeta), \phi_1(\zeta),
\phi_2(\zeta)$ and equating to zero coefficients at each powers of
$\zeta$ in (\ref{comp1}) we obtain (\ref{ex}).

{\bf Remark 1.} Let $\phi_1,...,\phi_{n+1}$ be linear independent
solutions of the linear system (\ref{lin}). We can construct a
system of the type (\ref{ex}) for each triplet of these solutions.
It is clear that all these systems are compatible. Therefore, our
systems have a lot of infinitesimal symmetries of hydrodynamic type.

Let us describe our pseudopotential written in the form
(\ref{pseudo}). One can derive differential equations for the
functions $f,g$ from (\ref{pspar}), (\ref{pspar1}) and
(\ref{pspar2}).

Define a function $q(\xi,u_1,...,u_n)$ as a solution of the
following system of PDEs
\begin{equation}   \label{par}
q_{\xi}=\frac{q^{s_1}(q-1)^{s_2}(q-u_1)^{s_3}...(q-u_n)^{s_{n+2}}}{\phi(q)},
\end{equation}
$$q_{u_i}=\frac{\phi(u_i)}{\phi(q)}\cdot\frac{q(q-1)(q-u_1)...\hat{i}...(q-u_n)}{u_i(u_i-1)(u_i-u_1)...\hat{i}...(u_i-u_n)},~~~ i=1,...,n.$$
The system (\ref{par}) is in involution iff the polynomial $\phi$
satisfies the linear system (\ref{lin}).

Let $\phi$, $\phi_1$ and $\phi_2$ be three linear independent
solutions of the system (\ref{lin}). Define functions
$f(\xi,u_1,...,u_n)$ and $g(\xi,u_1,...,u_n)$ as a solution of the
following system of PDEs
$$f_{\xi}=\frac{\phi_1(q)}{\phi(q)},\,   ~~~~
g_{\xi}=\frac{\phi_2(q)}{\phi(q)},$$
\begin{equation}   \label{ps}
f_{u_i}=-\frac{\mu_i(q)}{\phi(q)}q^{1-s_1}(q-1)^{1-s_2}(q-u_1)^{1-s_3}...(q-u_n)^{1-s_{n+2}},
\end{equation}
$$g_{u_i}=\frac{\nu_i(q)}{\phi(q)}q^{1-s_1}(q-1)^{1-s_2}(q-u_1)^{1-s_3}...(q-u_n)^{1-s_{n+2}},~~~i=1,..,n.$$
Here $\mu_i,\nu_i$ are defined by (\ref{mu}). It can be checked
straightforwardly that the system (\ref{ps}) is in involution.

{\bf Proposition 3.} If the functions $f,g$ are defined by
(\ref{ps}), then the system (\ref{pseudo}) is a pseudopotential for
the system (\ref{ex}).

{\bf Proof.} The formulas (\ref{ps}) imply
\begin{equation}
\label{ps1} f_{\xi} g_{u_i} - g_{\xi}
f_{u_i}=-\frac{\vartheta_i(q)}{\phi(q)}q^{1-s_1}(q-1)^{1-s_2}(q-u_1)^{1-s_3}...(q-u_n)^{1-s_{n+2}}
\end{equation}
where $\vartheta_i$ is defined by (\ref{mu}). Therefore, the space
of functions $\{f_{u_i},g_{u_i},f_{\xi} g_{u_i} - g_{\xi}
f_{u_i};i=1,...,n\}$ in the variable $\xi$ is isomorphic to the
space of polynomials in $q$ of degree less or equal to $n-1$. This
space is $n$-dimensional and we can apply Lemma 1. Let
$\phi(q)=\alpha_0+\alpha_1q+...+\alpha_nq^n$,
$\phi_1(q)=\beta_0+\beta_1q+...+\beta_nq^n$ and
$\phi_2(q)=\gamma_0+\gamma_1q+...+\gamma_nq^n$. Substituting
(\ref{ps}), (\ref{ps1}) into (\ref{comp}) we obtain (\ref{comp1}).
Calculating polynomials $\nu_i(q),\mu_i(q),\vartheta_i(q)$ in terms
of coefficients of polynomials $\phi(q), \phi_1(q), \phi_2(q)$ and
equating to zero coefficients at each powers of $q$ in (\ref{comp1})
we obtain (\ref{ex}).

{\bf Remark 2.} Let $s_3=...=s_{n+2}=1$ and
$\phi(q)=(q-u_1)...(q-u_n)$. Then the system (\ref{ex}) coincides
with a system from the paper \cite{odsok} (see Example 4 in
\cite{odsok}).

\section{The case $n=2$}

In this case each system (\ref{gener}) can be written in the form
\begin{equation}
\left(
\begin{array}{c}
v_t \\
\ \\
w_t
\end{array}
\right)+ \left(
\begin{array}{cc}
a & 0 \\
\ \\
0 & b
\end{array}
\right) \left(
\begin{array}{c}
v_x \\
\ \\
w_x
\end{array}
\right)+ \left(
\begin{array}{cc}
p&q \\
\ \\
r&s
\end{array}
\right) \left(
\begin{array}{c}
v_y \\
\ \\
w_y
\end{array}
\right)=0 \label{sys}
\end{equation}
in some suitable coordinates $v,w$. These systems were intensively
studied in the paper \cite{ferhus2}. In particular, it was proven in
this paper that the system (\ref{sys}) is integrable iff the
functions $a,b,p,q,r,s$ satisfy the following system of PDEs

\noindent {\bf Equations for $a$:}
\begin{eqnarray}
a_{vv}&=&\frac{\displaystyle{qa_vb_v+2qa_v^2+(s-p)a_va_w-ra_w^2}}{\displaystyle{(a-b)q}}+\frac{\displaystyle{a_vr_v}}{\displaystyle{r}}+\frac{\displaystyle{2a_vp_w-a_wp_v}}{\displaystyle{q}},
\nonumber \\
a_{vw}&=&a_v\frac{\displaystyle{a_w+b_w}}{\displaystyle{a-b}}+a_v(\frac{\displaystyle{q_w}}{q}+\frac{\displaystyle{r_w}}{r}),
\label{a} \\
a_{ww}&=&\frac{\displaystyle{qa_vb_v+(s-p)a_vb_w+ra_w^2}}{\displaystyle{(a-b)r}}+\frac{\displaystyle{a_vs_w}}{r}+\frac{\displaystyle{a_wq_w}}{q}.
\nonumber
\end{eqnarray}

\noindent  {\bf Equations for $b$:}
\begin{eqnarray}
\nonumber \\
b_{vv}&=&\frac{\displaystyle{ra_wb_w+(p-s)a_vb_w+qb_v^2}}{\displaystyle{(b-a)q}}+\frac{\displaystyle{b_wp_v}}{q}+\frac{\displaystyle{b_vr_v}}{r};
\nonumber \\
b_{vw}&=&b_w\frac{\displaystyle{a_v+b_v}}{\displaystyle{b-a}}+b_w(\frac{\displaystyle{q_v}}{q}+\frac{\displaystyle{r_v}}{r}),
\label{b} \\
b_{ww}&=&\frac{\displaystyle{ra_wb_w+2rb_w^2+(p-s)b_vb_w-qb_v^2}}{\displaystyle{(b-a)r}}+\frac{\displaystyle{b_wq_w}}{\displaystyle{q}}+\frac{\displaystyle{2b_ws_v-b_vs_w}}{\displaystyle{r}}.
\nonumber
\end{eqnarray}

\noindent  {\bf Equations for $p$:}
\begin{eqnarray}
p_{vv}&=&2\frac{r(a_vb_w-a_wb_v)+(s-p)a_vb_v}{(a-b)^2}+\frac{r_vp_v}{r}+\frac{p_vp_w}{q}+
\nonumber \\
&&\frac{\frac{r}{q}
(2q_va_w-2a_vq_w+a_wp_w)-b_vp_v+2r_va_w-2a_v(s_v+p_v+r_w)+\frac{p-s}{q}(2p_va_w-a_vp_w)}{b-a}\nonumber
\\
\ \nonumber \\
p_{vw}&=&2(s-p)\frac{a_vb_w}{(a-b)^2}-\frac{b_wp_v+(2s_w+p_w)a_v}{b-a}+{p_v}\left
(\frac{q_w}{q}+\frac{r_w}{r}\right ), \label{p} \\
\  \nonumber \\
p_{ww}&=&2\frac{q(a_wb_v-a_vb_w)+(s-p)a_wb_w}{(a-b)^2}+\frac{(p-s)b_wp_v-qb_vp_v-2rs_wa_w-ra_wp_w}{(b-a)r}+
\nonumber \\
&& \frac{p_vs_w}{r}+\frac{q_wp_w}{q}. \nonumber
\end{eqnarray}

\noindent  {\bf Equations for $s$:}
\begin{eqnarray}
s_{vv}&=&2\frac{r(a_wb_v-a_vb_w)+(p-s)a_vb_v}{(a-b)^2}+\frac{(s-p)a_vs_w-ra_ws_w-2qp_vb_v-qb_vs_v}{(a-b)q}+
\nonumber \\
&&\frac{p_vs_w}{q}+\frac{r_vs_v}{r}, \nonumber \\
\  \nonumber \\
s_{vw}&=&2(p-s)\frac{a_vb_w}{(a-b)^2}-\frac{a_vs_w+(2p_v+s_v)b_w}{a-b}+{s_w}\left
(\frac{q_v}{q}+\frac{r_v}{r}\right ), \label{s} \\
\  \nonumber \\
s_{ww}&=&2\frac{q(a_vb_w-a_wb_v)+(p-s)a_wb_w}{(a-b)^2}+\frac{q_ws_w}{q}+\frac{s_vs_w}{r}
+ \nonumber \\
&&\frac{\frac{q}{r}
(2r_wb_v-2b_wr_v+b_vs_v)-a_ws_w+2q_wb_v-2b_w(p_w+s_w+q_v)+\frac{s-p}{r}(2s_wb_v-b_ws_v)}{a-b}.
\nonumber
\end{eqnarray}

\noindent {\bf Equations for  $q$ and $r$:}
\begin{eqnarray}
qr_{ww}+rq_{ww}&=&2(p-s)\frac{(p-s)a_wb_w+q(a_vb_w-a_wb_v)}{(a-b)^2}
+q\frac{r_v}{r}\frac{qb_v+(s-p)b_w}{a-b}+   \nonumber \\
&&(s-p)\frac{2a_ws_w+2b_wp_w+b_wq_v}{a-b}+r\frac{(a_w-2b_w)q_w}{a-b}+
\nonumber \\
&&q\frac{a_wr_w+b_v(2p_w+2s_w+q_v)-2b_w(r_w+p_v+s_v)}{a-b}+ \nonumber \\
&&\frac{r}{q}q_w^2+\frac{q}{r}s_wr_v-q_wr_w+s_w(2p_w+q_v), \nonumber \\
\ \nonumber \\
q_{vw}&=&(s-p)\frac{qa_vb_v+(s-p)a_vb_w+ra_wb_w}{r(a-b)^2}
+\frac{q_vq_w}{q}+\frac{p_vs_w}{r}+                     \nonumber \\
&&
\frac{a_v(rq_w+qr_w)+(s-p)(a_vs_w+b_wp_v)+ra_ws_w+qp_vb_v}{r(a-b)},
\nonumber \\
\  \label{qr} \\
r_{vw}&=&(p-s)\frac{ra_wb_w+(p-s)a_vb_w+qa_vb_v}{q(a-b)^2}
+\frac{r_vr_w}{r}+\frac{p_vs_w}{q}+                     \nonumber \\
&&
\frac{b_w(rq_v+qr_v)+(p-s)(a_vs_w+b_wp_v)+ra_ws_w+qp_vb_v}{q(b-a)},
\nonumber \\
\ \nonumber \\
qr_{vv}+rq_{vv}&=&2(s-p)\frac{(s-p)a_vb_v+r(a_vb_w-a_wb_v)}{(a-b)^2}
+r\frac{q_w}{q}\frac{ra_w+(p-s)a_v}{b-a}+   \nonumber \\
&&(p-s)\frac{2b_vp_v+2a_vs_v+a_vr_w}{b-a}+q\frac{(b_v-2a_v)r_v}{b-a}+
\nonumber \\
&&r\frac{b_vq_v+a_w(2s_v+2p_v+r_w)-2a_v(q_v+s_w+p_w)}{b-a}+ \nonumber \\
&&\frac{q}{r}r_v^2+\frac{r}{q}p_vq_w-r_vq_v+p_v(2s_v+r_w); \nonumber
\end{eqnarray}
Suppose that the system (\ref{sys}) possesses a pseudopotential of
the form
\begin{equation}
\psi_t=f(\psi_y, \ v, \ w), ~~~ \psi_x=g(\psi_y, \ v, \ w).
\label{Lax}
\end{equation}
Another remarkable result of the paper \cite{ferhus2} is the
following

{\bf Theorem.} The class of two-component (2+1)-dimensional systems
of hydrodynamic type possessing infinitely many hydrodynamic
reductions coincides with the class of systems possessing a
pseudopotential of the form (\ref{Lax}).

Recall the prove of this theorem (see \cite{ferhus2} for details).
Writing out the consistency condition $\psi_{tx}=\psi_{xt}$,
expressing $v_t, \ w_t$  by virtue of (\ref{sys}) and equating to
zero coefficients at $v_x, v_y, w_x, w_y$, one arrives at the
following expressions for the first derivatives $f_v, f_w, f_{\xi}$
and $g_{\xi}$ (we adopt the notation $\xi \equiv \psi_y$):
\begin{equation}
\begin{array}{c}
f_v=-a\ g_v, ~~~ f_w=-b\ g_w, \\
\ \\
f_{\xi}=\frac{{b\left(p+r\frac{g_w}{g_v}\right)-a\left(s+q\frac{g_v}{g_w}\right)}}{{a-b}},
\end{array}
\label{f}
\end{equation}
and
\begin{equation}
  g_{\xi}=\frac{{s+q\frac{g_v}{g_w}-p-r\frac{g_w}{g_v}}}{{a-b}}.
\label{g}
\end{equation}
The consistency conditions of the equations (\ref{f}) imply the
following expressions for the second partial derivatives $g_{vw},
g_{vv}, g_{ww}$:
\begin{equation}
\begin{array}{c}
g_{vw}=\frac{a_w}{b-a}\ g_v+\frac{b_v}{a-b}\ g_w, \\
\ \\
g_{vv}=\frac{g_v[g_w^2(r(b_v-a_v)+(a-b)r_v)+g_vg_w((a-b)p_v+(s-p)a_v-ra_w)+qa_vg_v^2]}{(a-b)
r g_w^2}, \\
\  \\
g_{ww}=\frac{g_w[g_v^2(q(a_w-b_w)+(b-a)q_w)+g_vg_w((b-a)s_w+(p-s)b_w-qb_v)+rb_wg_w^2]}{(b-a)
q g_v^2}.
\end{array}
\label{g2}
\end{equation}
The compatibility conditions of the equations (\ref{g}), (\ref{g2})
for $g$, namely, the conditions  $g_{\xi vv}=g_{vv \xi}, \ g_{\xi
vw}=g_{vw \xi} $, etc., are of the form $P(g_v, g_w)=0$, where $P$
denotes a rational expression in $g_v, \ g_w$ whose coefficients are
functions of $a, b, p, q, r, s$ and their partial derivatives  up to
the second order. Equating all these expressions to zero (they are
required to be zero identically in $g_v, g_w$), one obtains the  set
of conditions which are necessary and sufficient for the existence
of a pseudopotential of the form (\ref{Lax}). It turns out that
these conditions identically coincide with the integrability
conditions (\ref{a}) - (\ref{qr}). Thus, any system satisfying the
integrability conditions  (\ref{a}) - (\ref{qr}) possesses a
pseudopotential of the form (\ref{Lax}).

Based on these results of \cite{ferhus2} we can prove the following

{\bf Proposition 4.} If a system (\ref{sys}) corresponds to a
general solution of (\ref{a}) - (\ref{qr}), then in suitable
coordinates it is equivalent to some system of the form (\ref{ex})
(in the case $n=2$) and possesses a pseudopotential of the form
(\ref{ps}).

{\bf Proof.} Let us study the system (\ref{f})-(\ref{g2}). We denote
$q=\frac{g_w}{g_v}$. Calculation shows that
$$\frac{q_{\xi\xi}}{q_{\xi}^2}=-\frac{1}{q}+\frac{P_4(q)}{P_5(q)}$$
where $P_4$, $P_5$ are polynomials in $q$ of degree 4 and 5
correspondingly. Moreover, $P_5$ has distinct roots for general
solution of the system (\ref{a}) - (\ref{qr}). Write
$$\frac{q_{\xi\xi}}{q_{\xi}^2}=-\frac{1}{q}+\sum_{1\leq
i\leq5}\frac{s_i}{q-\lambda_i}$$ where $s_i$, $\lambda_i$ are some
functions in $v,w$. Calculation shows that $\sum_{1\leq
i\leq5}s_i=\lim_{q\to\infty}\frac{qP_4(q)}{P_5(q)}=4$. After
integration we get
$$q_{\xi}=\frac{C}{q}\prod_{1\leq
i\leq5}(q-\lambda_i)^{s_i}.$$ Note also that
$$f_{\xi}=\frac{P_2(q)}{q},~~~ g_{\xi}=\frac{Q_2(q)}{q}$$ where
$P_2$, $Q_2$ are quadratic polynomials in $q$ and
$$q_v=\frac{S_2(q)}{q},~~~ q_w=qG_2(q)$$ for some quadratic polynomials $S_2$, $G_2$.
Let us make a substitution of the form
\begin{equation}   \label{subs}
q\to\frac{\alpha q+\beta}{\gamma q+\delta}
\end{equation}
where
$\alpha,\beta,\gamma,\delta$ are functions in $v,w$. After that we
obtain
\begin{equation}   \label{t}
q_{\xi}=\frac{1}{\phi(q)}\prod_{1\leq i\leq5}(q-\rho_i)^{s_i},$$
$$q_v=\frac{S_4(q)}{\phi(q)}, q_w=\frac{G_4(q)}{\phi(q)}
\end{equation}
$$f_{\xi}=\frac{\phi_1(q)}{\phi(q)},
g_{\xi}=\frac{\phi_2(q)}{\phi(q)}$$ where $\phi(q)$, $\phi_1(q)$,
$\phi_2(q)$ are quadratic polynomials in $q$ and $S_4$, $G_4$ are
polynomials of degree 4.

Note that after appropriate substitution (\ref{subs}) we can assume
$\rho_1=0$, $\rho_2=1$ and $\rho_5=\infty$. Moreover, after change
of variables $v=v(u_1,u_2)$, $w=w(u_1,u_2)$ we can assume
$\rho_3=u_1$, $\rho_4=u_2$. Compatibility conditions for the system
(\ref{t}) with these assumptions imply that this system has a form
(\ref{par}) with $n=2$ where $\phi$ satisfies (\ref{lin}). In
particular, $s_1,...,s_5$ must be constant. Moreover, coefficients
of $\phi$, $\phi_1$, $\phi_2$ must satisfy the same system of linear
equations. Indeed, we can swap $\phi,\phi_1,\phi_2$ by exchanging
the role of $t,x,y$ in (\ref{Lax}). After that it is easy to find
$f_{u_i},g_{u_i}$ in the form (\ref{ps}).

\section{Elliptic case}

Fix a complex parameter $\tau$ such that Im$\tau>0$. Define a
function $\theta(z)$ by
$$\theta(z)=\sum_{m\in\Z}(-1)^me^{-2\pi
i(mz+\frac{m(m-1)}{2}\tau)}.$$ It is clear that $\theta(z)$ is the
entire function satisfying the following relations:
$$\theta(z+1)=\theta(z),~~~\theta(z+\tau)=-e^{-2\pi iz}\theta(z).$$
Moreover, each entire function in one variable satisfying these
relations is proportional to $\theta(z)$. We have also
$\theta(-z)=-e^{-2\pi iz}\theta(z)$ and the only zero of the
function $\theta(z)$ modulo 1 and $\tau$ is $z=0$.

Let $\Theta_{n,c}(\tau)$ be the space of the entire functions in one
variable satisfying the following relations:
$$f(z+1)=f(z),~~~f(z+\tau)=(-1)^ne^{-2\pi i(nz-c)}f(z).$$
Here $n\in\N$ and $c\in\C$. It is known that $\dim
\Theta_{n,c}(\tau)=n$, every function $f\in \Theta_{n,c}(\tau)$ has
exactly $n$ zeros modulo 1 and $\tau$ (counting according to their
multiplicities), and the sum of these zeros is equal to $c$ modulo 1
and $\tau$. We have $\theta(z)\in \Theta_{1,0}(\tau)$.

Define a function $F(\zeta, u_1,\dots, u_n)$ as a solution of the
following systems of PDEs
$$F_{\zeta}=\frac{\phi(\zeta)}{\theta(\zeta-u_1)...\theta(\zeta-u_n)},$$
\begin{equation}   \label{psparell}
F_{u_i}=-\frac{1}
{\theta(u_i-u_1)...\hat{i}...\theta(u_i-u_n)}\cdot\frac{\theta(\zeta-u_i+\eta)}{\theta(\eta)\theta(\zeta-u_i)}\cdot\phi(u_i),~~~i=1,...,n.
\end{equation}
Here $\eta$ is a constant and
$\phi\in\Theta_{n,u_1+...+u_n-\eta}(\tau)$ as a function in $\zeta$.
This means that $\phi(\zeta)$ is the entire function in $\zeta$ and
\begin{equation}
\label{ell}\phi(\zeta+1)=\phi(\zeta),~~~\phi(\zeta+\tau)=(-1)^ne^{-2\pi
i(n\zeta-u_1-...-u_n+\eta)}\phi(\zeta).\end{equation} We assume that
$\eta$ is nonzero modulo 1 and $\tau$. The system (\ref{psparell})
is in involution iff the function $\phi$ satisfies the following
system of PDEs
\begin{equation}   \label{linell}
\phi_{u_i}(\zeta)=\phi(u_i)\frac{\theta(\zeta-u_1)...\hat{i}...\theta(\zeta-u_n)}{\theta(u_i-u_1)...\hat{i}...\theta(u_i-u_n)}\cdot
\frac{\theta(\zeta-u_i+\eta)}{\theta(\eta)}\times
\end{equation}
$$\left(\frac{\theta^{\prime}(\zeta-u_i)}{\theta(\zeta-u_i)}-\frac{\theta^{\prime}(\zeta-u_i+\eta)}{\theta(\zeta-u_i+\eta)}\right)-
\frac{\theta^{\prime}(\zeta-u_i)}{\theta(\zeta-u_i)}\phi(\zeta),~~~
i=1,...,n.$$ It is clear that if
$\phi\in\Theta_{n,u_1+...+u_n-\eta}(\tau)$, then the equations
(\ref{linell}) are compatible with (\ref{ell}) and the right hand
side of (\ref{linell}) is an entire function in $\zeta$. Therefore,
(\ref{linell}) is a well-defined system of linear PDEs for
coefficients of $\phi$ with respect to some basis in the
$n$-dimensional space $\Theta_{n,u_1+...+u_n-\eta}(\tau)$. It can be
checked straightforwardly that this system is in involution.
Therefore, there are $n$ linear independent solutions. Moreover,
this system can be solved explicitly. Namely, let
\begin{equation}   \label{sol}
\psi_i(\zeta)=\theta(\zeta-u_1)...\hat{i}...\theta(\zeta-u_n)\cdot\theta(\zeta-u_i+\eta),~~~i=1,...,n.
\end{equation}
It can be checked straightforwardly that these functions are linear
independent solutions of the system (\ref{linell}).

Assume that $n\geq3$. Let $\phi$, $\phi_1$ and $\phi_2$ be three
linear independent solutions of the system (\ref{linell}). We assume
that $\phi,\phi_1,\phi_2\in\Theta_{n,u_1+...+u_n-\eta}(\tau)$ as
functions in $\zeta$. Define functions $G(\zeta, u_1,\dots, u_n)$
and $H(\zeta, u_1,\dots, u_n)$ similarly to (\ref{psparell}) by
$$G_{\zeta}=\frac{\phi_1(\zeta)}{\theta(\zeta-u_1)...\theta(\zeta-u_n)},$$
\begin{equation}   \label{psparell1}
G_{u_i}=-\frac{1}
{\theta(u_i-u_1)...\hat{i}...\theta(u_i-u_n)}\cdot\frac{\theta(\zeta-u_i+\eta)}{\theta(\eta)\theta(\zeta-u_i)}\cdot\phi_1(u_i),~~~i=1,...,n
\end{equation}
for the function $G$ and
$$H_{\zeta}=\frac{\phi_2(\zeta)}{\theta(\zeta-u_1)...\theta(\zeta-u_n)},$$
\begin{equation}   \label{psparell2}
H_{u_i}=-\frac{1}
{\theta(u_i-u_1)...\hat{i}...\theta(u_i-u_n)}\cdot\frac{\theta(\zeta-u_i+\eta)}{\theta(\eta)\theta(\zeta-u_i)}\cdot\phi_2(u_i),~~~i=1,...,n
\end{equation}
for the function $H$.

{\bf Proposition 5.} If the functions $F,G,H$ are defined by
(\ref{psparell}), (\ref{psparell1}) and (\ref{psparell2}), then the
system (\ref{pseudopar}) defines a pseudopotential for some system
of the form (\ref{genern}) with $m=n$.

{\bf Proof.} Equations (\ref{psparell}), (\ref{psparell1}),
(\ref{psparell2}) imply
$$H_{\zeta}G_{u_i}-G_{\zeta}H_{u_i}=\frac{\vartheta_i(\zeta)}{\theta(\zeta-u_1)...\theta(\zeta-u_n)},$$
\begin{equation}   \label{exp}
F_{\zeta}H_{u_i}-H_{\zeta}F_{u_i}=\frac{\nu_i(\zeta)}{\theta(\zeta-u_1)...\theta(\zeta-u_n)},
\end{equation}
$$G_{\zeta}F_{u_i}-F_{\zeta}G_{u_i}=\frac{\mu_i(\zeta)}{\theta(\zeta-u_1)...\theta(\zeta-u_n)}$$
where functions $\vartheta_i(\zeta),\nu(\zeta),\mu(\zeta)$ are
defined by
$$\vartheta_i(\zeta)=\frac{1}{\theta(u_i-u_1)...\hat{i}...\theta(u_i-u_n)}\cdot\frac{\theta(\zeta-u_i+\eta)}{\theta(\eta)\theta(\zeta-u_i)}\cdot
(\phi_2(u_i)\phi_1(\zeta)-\phi_1(u_i)\phi_2(\zeta)),$$
\begin{equation}   \label{muell}
\nu_i(\zeta)=\frac{1}{\theta(u_i-u_1)...\hat{i}...\theta(u_i-u_n)}\cdot\frac{\theta(\zeta-u_i+\eta)}{\theta(\eta)\theta(\zeta-u_i)}\cdot
(\phi(u_i)\phi_2(\zeta)-\phi_2(u_i)\phi(\zeta)),
\end{equation}
$$\mu_i(\zeta)=\frac{1}{\theta(u_i-u_1)...\hat{i}...\theta(u_i-u_n)}\cdot\frac{\theta(\zeta-u_i+\eta)}{\theta(\eta)\theta(\zeta-u_i)}\cdot
(\phi_1(u_i)\phi(\zeta)-\phi(u_i)\phi_1(\zeta)).$$ It is clear that
$\vartheta_i(\zeta)$, $\nu_i(\zeta)$,
$\mu_i(\zeta)\in\Theta_{n,u_1+...+u_n-2\eta}(\tau)$. Therefore, the
space of functions $\{H_{\zeta}G_{u_i}-G_{\zeta}H_{u_i},
F_{\zeta}H_{u_i}-H_{\zeta}F_{u_i},
G_{\zeta}F_{u_i}-F_{\zeta}G_{u_i}; i=1,...,n\}$ in the variable
$\zeta$ is isomorphic to the space
$\Theta_{n,u_1+...+u_n-2\eta}(\tau)$. This space is $n$-dimensional
and we can apply Lemma 2.

Let us construct the system possessing pseudopotential defined by
(\ref{psparell}), (\ref{psparell1}), (\ref{psparell2}) explicitly.

{\bf Proposition 6.} Let
$$\phi(\zeta)=\sum_{i=1}^n\alpha_i\psi_i(\zeta),~~~\phi_1(\zeta)=\sum_{i=1}^n\beta_i\psi_i(\zeta),~~~\phi_2(\zeta)=\sum_{i=1}^n\gamma_i\psi_i(\zeta)$$
where $\psi_i(\zeta)$ are given by (\ref{sol}) and
$\alpha_i,\beta_i,\gamma_i$ are constants. Then the system with
pseudopotential defined by (\ref{psparell}), (\ref{psparell1}),
(\ref{psparell2}) can be written in the form
\begin{equation}   \label{exell}
\sum_{i}\frac{\theta(u_j-u_i+\eta)}{\theta(u_j-u_i)}((\gamma_i\alpha_j-\gamma_j\alpha_i)(u_{it}-u_{jt})+
(\alpha_i\beta_j-\alpha_j\beta_i)(u_{ix}-u_{jx})+(\beta_i\gamma_j-\beta_j\gamma_i)(u_{iy}-u_{jy}))=0.
\end{equation}
Here summation is made by $i$ subject to the constrains $1\leq i\leq
n$, $i\ne j$ with fixed $j$. For each $j=1,...,n$ we have an
equation.

{\bf Proof.} Substituting (\ref{exp}) into (\ref{comppar}) we obtain
\begin{equation}   \label{comp1ell}
 \sum_{i=1}^n \nu_i(\zeta)u_{it}+\sum_{i=1}^n
\mu_i(\zeta)u_{ix}+\sum_{i=1}^n \vartheta_i(\zeta)u_{iy}=0.
\end{equation}
Calculating functions $\nu_i(\zeta),\mu_i(\zeta),\vartheta_i(\zeta)$
in terms of functions $\phi(\zeta), \phi_1(\zeta), \phi_2(\zeta)$
and evaluating (\ref{comp1ell}) at $\zeta=u_j$ we obtain
(\ref{exell}).

{\bf Remark 3.} If $n=3$, then the system (\ref{exell}) reduces to
the trivial system
$$u_{1t}=u_{3t},~~~u_{2x}=u_{1x},~~~u_{2y}=u_{3y}.$$
Therefore, one can assume $n>3$.

{\bf Remark 4.} The system (\ref{exell}) is invariant under
translation $u_i\to u_i+v$ for an arbitrary function $v$. Therefore,
one can reduce the number of unknown functions in the system setting
$u_n=0$.

{\bf Remark 5.} The system (\ref{exell}) is written in the form
(\ref{genern}) and can not be written in the form (\ref{gener}). In
particular, it does not belong to the class of the systems studied
in the paper \cite{ferhus3}.

Let us describe our pseudopotential written in the form
(\ref{pseudo}). One can derive differential equations for the
functions $f,g$ from (\ref{psparell}), (\ref{psparell1}) and
(\ref{psparell2}).

Define a function $q(\xi,u_1,...,u_n)$ as a solution of the
following system of PDEs
\begin{equation}   \label{parell}
q_{\xi}=\frac{\theta(q-u_1)...\theta(q-u_n)}{\phi(q)},
\end{equation}
$$q_{u_i}=\frac{\phi(u_i)}{\phi(q)}\cdot\frac{\theta(q-u_1)...\hat{i}...\theta(q-u_n)}{\theta(u_i-u_1)...\hat{i}...\theta(u_i-u_n)}\cdot
\frac{\theta(q-u_i+\eta)}{\theta(\eta)},~~~ i=1,...,n.$$ The system
(\ref{parell}) is in involution iff the function
$\phi\in\Theta_{n,u_1+...+u_n-\eta}(\tau)$ satisfies (\ref{linell}).

Let $\phi$, $\phi_1$, $\phi_2\in\Theta_{n,u_1+...+u_n-\eta}(\tau)$
be three linear independent solutions of the system (\ref{linell}).
Define functions $f(\xi,u_1,...,u_n)$ and $g(\xi,u_1,...,u_n)$ as a
solution of the following system of PDEs
$$f_{\xi}=\frac{\phi_1(q)}{\phi(q)},\,   ~~~~
g_{\xi}=\frac{\phi_2(q)}{\phi(q)},$$
\begin{equation}   \label{psell}
f_{u_i}=-\frac{\mu_i(q)}{\phi(q)},~~~g_{u_i}=\frac{\nu_i(q)}{\phi(q)},~~~i=1,..,n.
\end{equation}
Here $\mu_i,\nu_i$ are defined by (\ref{muell}). It can be checked
straightforwardly that the system (\ref{psell}) is in involution.

{\bf Proposition 7.} If the functions $f,g$ are defined by
(\ref{psell}), then the system (\ref{pseudo}) is a pseudopotential
for some system of the form (\ref{genern}).

{\bf Proof.} The formulas (\ref{psell}) imply
\begin{equation}
\label{ps1ell} f_{\xi} g_{u_i} - g_{\xi}
f_{u_i}=-\frac{\vartheta_i(q)}{\phi(q)}
\end{equation}
where $\vartheta_i$ is defined by (\ref{muell}). Therefore, the
space of functions $\{f_{u_i},g_{u_i},f_{\xi} g_{u_i} - g_{\xi}
f_{u_i};i=1,...,n\}$ in the variable $\xi$ is isomorphic to the
space $\Theta_{n,u_1+...+u_n-2\eta}(\tau)$ in the variable $q$. This
space is $n$-dimensional and we can apply Lemma 2.

\section{Discussion}

In this paper we suggest a construction of $n$-component
(2+1)-dimensional hydrodynamic type systems possessing
pseudopotential. Let us outline some features of this construction.
It is clear from (\ref{psecon1}), (\ref{psecon2}) that there exists
a polynomial $S(u,v)$ of degree $n$ such that
$S(f_{\xi},g_{\xi})=0$. If a curve $K=\{(u,v); S(u,v)=0\}$ is
rational, then it is natural to write $f_{\xi},g_{\xi}$ in the form
(\ref{ps}) where $\phi,\phi_1,\phi_2$ are polynomials of degree $n$.
It turns out that the coefficients of each polynomial satisfy the
same system of linear PDEs. A similar construction exists if the
curve $K$ is elliptic, but in this case we get an overdetermined
integrable system with $n$ equations for $n-1$ unknowns. The natural
conjecture is that there exists a similar construction in the case
$g>1$ where $g$ is the genus of $K$ and the number of equations in
the corresponding integrable system should be $g$ plus the number of
unknowns. It would be interesting to verify this conjecture and find
these systems explicitly as well as their degenerations. Some (or
may be even all in the case $g>0$) of these systems could be
reductions of the universal Whitham hierarchy \cite{kr4}.

Our construction provides a general solution to the classification
problem of two-component integrable hydrodynamic type systems. The
full classification of these systems as well as their detailed study
will be the subject of a separate paper.

\vskip.3cm \noindent {\bf Acknowledgments.} I am grateful to
E.Ferapontov and V.Sokolov for useful discussions. E.Ferapontov
attracted my attention to his paper \cite{ferhus2} which became a
base of this research. Working with V.Sokolov at our paper
\cite{odsok} I understood the important role of pseudopotentials in
the theory of hydrodynamic type systems.

\end{document}